\documentclass[12pt]{article}
\usepackage{amssymb,amsmath,latexsym}

\def\proof{{\it Proof}.\ }

\def\wbull{\hfill\vrule height .9ex width .8ex depth -.1ex}

\newtheorem{formula}{}[section]
\newtheorem{proposition}[formula]{Proposition}
\newtheorem{definition}[formula]{Definition}
\newtheorem{corollary}[formula]{Corollary}

\newtheorem{lemma}[formula]{Lemma}
\newtheorem{theorem}[formula]{Theorem}

\def\thrm{\begin{theorem}}
\def\thrml#1{\begin{theorem}\label{#1}}
\def\ethrm{\end{theorem}}
\def\prp{\begin{proposition}}
\def\prp#1{\begin{proposition}\label{#1}}
\def\eprp{\end{proposition}}
\def\dfntn{\begin{definition}}
\def\dfntnl#1{\begin{definition}\label{#1}}
\def\edfntn{\end{definition}}
\def\nmrt{\begin{enumerate}}
\def\enmrt{\end{enumerate}}

\def\qtn{\begin{equation}}
\def\qtnl#1{\begin{equation}\label{#1}}
\def\eqtn{\end{equation}}
\def\lmm{\begin{lemma}}
\def\lmml#1{\begin{lemma}\label{#1}}
\def\elmm{\end{lemma}}
\def\crllr{\begin{corollary}}
\def\crllrl#1{\begin{corollary}\label{#1}}
\def\ecrllr{\end{corollary}}
\def\css{\begin{cases}}
\def\ecss{\end{cases}}

\title{ \bf{There exists no distance-regular graph with intersection array 
$\{56,36,9;1,3,48\}$}}
\author{
Alexander Gavrilyuk
\thanks{
Partially supported by the Russian Foundation for Basic Research (project no. 08-01-00009).}\\[-1pt]
\small Ural Division of the Russian Academy of Sciences\\[-3pt]
\small Institute of Mathematics and Mechanics\\[-3pt]
\small ul.~S.~Kovalevskoi~16, Yekaterinburg, 620990 Russia\\[-3pt]
{\tt \small alexander.gavriliouk@gmail.com}\\[-3pt]}

\begin{document}

\maketitle

\begin{abstract}
We prove that a distance-regular graph with intersection array 
$\{56,36,9;1,3,48\}$ does not exist. This intersection array
is from the table of feasible parameters for distance-regular graphs
in "Distance-regular graphs"\ by A.E. Brouwer, A.M. Cohen, A. Neumaier.
\medskip

\end{abstract}
\newpage

\section{Introduction}

The purpose of this paper is to prove the next

\thrml{T}
The array $\{56,36,9;1,3,48\}$ cannot be realized as 
the intersection array of a distance-regular graph.
\ethrm

This array is from the table of feasible parameters for distance-regular graphs
in "Distance-regular graphs"\ by A.E. Brouwer, A.M. Cohen, A. Neumaier 
(see~\cite[p. 429]{BCN}). 

The important tool of our theorem's proof is the Koolen-Park inequality 
(see next section). This inequality shows the largest coclique 
in the vertex neighborhood of hypothetic distance-regular graph $\Gamma$ 
with intersection array from Theorem \ref{T} has size 3 
(i.e., $\Gamma$ does not contain a 4-claw). 
Using this observation, the possible neighborhoods of vertices of $\Gamma$
are determined. For each of them, we will get a contradiction.

\section{Definitions and preliminaries}

We consider only finite undirected graphs without loops or multiple edges.
Let $\Gamma$ be a connected graph. The \emph{distance} ${\rm d}(u,w)$ between any
two vertices $u$ and $w$ of $\Gamma$ is the length of a shortest path from $u$ to
$w$ in $\Gamma$. The \emph{diameter} ${\rm diam}(\Gamma)$ of $\Gamma$ is the maximal
distance occurring in $\Gamma$.

For a subset $A$ of the vertex set of $\Gamma$, we will also write $A$ for
the subgraph of $\Gamma$ induced by $A$.
For a vertex $u$ of $\Gamma$, define $\Gamma_i(u)$ to be the set
of vertices that are at distance $i$ from $u$ ($0\le i\le {\rm diam}(\Gamma))$.
The subgraph $\Gamma_1(u)$ is called the \emph{neighborhood}) of a vertex $u$
(and it will be simply denoted as $\Gamma(u)$) 
and the \emph{degree} of $u$ is the number of neighbors of $u$, 
i.e., $|\Gamma(u)|$.
A graph is \emph{regular} with degree $k$
if the degree of each of its vertices is $k$.

For the vertices $u_1,u_2,\ldots,u_s\in \Gamma$, define 
$\Gamma(u_1,u_2,\ldots,u_s)$ be the set of vertices of $\cap_{i=1}^s \Gamma(u_i)$.
For two vertices $u,w\in \Gamma$ with ${\rm d}(u,w)=2$, 
the subgraph $\Gamma(u,w)$ is called the $\mu$-\emph{subgraph} 
of vertices $u,w$. 
\medskip

A connected graph $\Gamma$ with diameter $d={\rm diam}(\Gamma)$ is \emph{distance-regular}
if there are integers $b_i$, $c_i$ ($0\le i\le d$) such that, for any two vertices $u,w\in \Gamma$
with ${\rm d}(u, w)=i$, there are exactly $c_i$ neighbors of $w$ in $\Gamma_{i-1}(u)$ and $b_i$
neighbors of $w$ in $\Gamma_{i+1}(u)$
(we assume that $\Gamma_{-1}(u)$ and $\Gamma_{d+1}(u)$ are empty sets).
In particular, a distance-regular graph $\Gamma$ is
regular with degree $b_0$, $c_1=1$ and $c_2=\mu(\Gamma)$.
For each vertex $u\in \Gamma$ and $0\le i\le d$, the subgraph $\Gamma_i(u)$ is regular 
with degree $a_i=b_0-b_i-c_i$. The numbers $a_i$, $b_i$, $c_i$ ($0\le i\le d$) are called 
the \emph{intersection numbers} and 
the array $\{b_0,b_1,\ldots,b_{d-1};c_1,c_2,\ldots,c_{d}\},$ is called
the \emph{intersection array} of the distance-regular graph $\Gamma$.
\medskip

A $c$-\emph{clique} $C$ of $\Gamma$ is a complete subgraph 
(i.e., every two vertices of $C$ are adjacent) of $\Gamma$ 
with exactly $c$ vertices. 
We say that $C$ is a clique if it is a $c$-clique for certain $c$.
A coclique $C$ of $\Gamma$ is an induced subgraph of $\Gamma$ 
with empty edge set. 
We say a coclique is a $c$-\emph{coclique} if it has exactly $c$ vertices.
\medskip

The following lemma is due to J.H. Koolen and J. Park~\cite{KP} 
(see also~\cite{CDG}).

\lmml{KPineq}
Let $\Gamma$ be a distance-regular graph and, for a vertex $x\in \Gamma$,
the neighborhood of $x$ contains a coclique of size $c\ge 2$. 
Then
$$c_2-1\ge \displaystyle {{\frac {c (a_1+1)-b_0}{{c\choose 2}}}}.$$
\elmm

\section{A proof of theorem}

Let $\Gamma$ be a distance-regular graph with intersection array
$\{56,36,9;1,3,48\}$. The intersection number $a_1$ of $\Gamma$ 
equals $56-36-1=19$. Fix an arbitrary vertex $\infty$ of $\Gamma$ 
and denote the subgraph $\Gamma(\infty)$ by $\Delta$.
In particular, the graph $\Delta$ is regular of degree $19$ and, 
for each pair of nonadjacent vertices $x,y$ of $\Delta$, 
$|\Delta(x,y)|\le 2$ holds.

\lmml{KPcocl}
The largest coclique of $\Delta$ has size $3$.
Moreover, each vertex of $\Delta$ belongs to a maximal coclique 
of size $3$. 
\elmm
\proof 
Let $\Delta$ contain a $c$-coclique. 
It is easy to see that each vertex of $\Delta$ belongs to a coclique 
of size at least $b_0/(a_1+1)=56/20>2$. Hence, we may assume $c\ge 3$.
If $c=4$, then, by Lemma \ref{KPineq}, $3-1\ge (4(19+1)-56)/6=4$, a contradiction. 
\wbull
\medskip

Let the vertices $x_1,x_2,x_3\in \Delta$ induce a 3-coclique.
Denote the vertex set of ${\displaystyle \Delta(x_i)-\cup_{j\ne i}\Delta(x_j)}$ by $X_i$.

\lmml{Cases}
Without loss of generality, one of the following cases holds.

\noindent $(1)$ $\Delta(x_1,x_2)=\{u,w\}$, 
$\Delta(x_1,x_3)=\{p\}$ and $\Delta(x_2,x_3)=\{q\}$,

\noindent $(2)$ $\Delta(x_1,x_2)=\{u,w\}$ 
and $\Delta(x_2,x_3)=\{p,q\}$,

\noindent $(3)$ $\Delta(x_1,x_2,x_3)=\{u\}$, 
$\Delta(x_1,x_2)=\{u,p\}$, $\Delta(x_2,x_3)=\{u,q\}$,

\noindent $(4)$ $\Delta(x_1,x_2,x_3)=\{u,w\}$.
\elmm
\proof For each pair of distinct indices $i,j\in\{1,2,3\}$, we denote 
$|\Delta(x_i,x_j)-\Delta(x_1,x_2,x_3)|$ by $\delta_{ij}$ and 
$|\Delta(x_1,x_2,x_3)|$ by $\delta$. Then we have 
$|X_i|=19-\sum_{j,j\ne i}\delta_{ij}-\delta$ and 
$|\Delta|=3+\delta+\sum_{i<j}\delta_{ij}+\sum_{i} |X_i|$.
Hence, $60-2\delta-\sum_{i<j}\delta_{ij}=56$ 
and $2\delta+\sum_{i<j}\delta_{ij}=4$. 
If $\delta=0$, then either $\delta_{ij}=2$ for two pair of indices (and we have Case (2))
or $\delta_{ij}=1$ or 2 and we have Case (1). 
If $\delta=1$, then $\delta_{ij}=1$ for two pair of indices and we have Case (3).
If $\delta=2$, then $\delta_{ij}=0$ and we have Case (4).
The lemma is proved.
\wbull

\lmml{Claim12}
The following hold.

\noindent $(1)$ $X_i$ is a clique$;$

\noindent $(2)$ For a vertex $z$ of $\Delta(x_i,x_j)$, 
either $X_i\subset \Delta(z)$ or $|X_i\cap \Delta(z)|\le 1$ holds$.$
\elmm
\proof (1) If $X_i$ contains a 2-coclique $\{a,b\}$, 
then the vertex set of ${\displaystyle\cup_{j\ne i}\{x_j\}\cup \{a,b\}}$ 
induce a 4-coclique in $\Delta$. This contradicts Lemma \ref{KPcocl}.

(2) Suppose that $1<|X_i\cap \Delta(z)|<|X_i|$.
Then, for a vertex $y\in X_i-\Delta(z)$, 
the $\mu$-subgraph of $z,y$ contains the vertices $\infty,x_i$ 
and $|X_i\cap \Delta(z)|\ge 2$ vertices of $X_i$, which is impossible.
\wbull

\lmml{Claim3}
\noindent Let $\{a,b\}$ be an edge of $\Delta$ such that $|\Delta(a,b)|\le 1$, 
$c$ be a vertex of $\Gamma_2(\infty)\cap \Gamma(a,b)$ and 
$d$ be a vertex of $\Gamma(\infty,c)-\{a,b\}$.
Then $d\in\Delta(a)\cup \Delta(b)$ and 
if $d\not\sim b$, then $\Delta(b,d)=\{a\}$.
\elmm
\proof The subgraph $\Delta-\Delta(a)\cup \Delta(b)$ 
contains at most $56-2\cdot 17-2-1=19$ vertices.
Since ${\rm d}(\infty,c)=2$, $\Gamma(\infty,c)$ contains 
$a,b$ and one more vertex, say, $d$. If $d\not\sim a$ and $d\not\sim b$, 
then $|\Gamma_3(c)\cap \Gamma(\infty)|\le 56-(2\cdot 17+1+(19-2))<b_2=9$,
a contradiction.
If $d\in \Delta(a)-\Delta(b)$, then $\Gamma(b,d)=\{a,c,\infty\}$.
\wbull

\lmml{Case4}
Case $(4)$ is impossible.
\elmm
\proof Note that $|X_i|=17$ for $i=1,2,3$. 
By Lemma \ref{Claim12}(2), we may suppose $X_1\subset \Delta(u)$.
Since $|\Gamma(u,w)|\ge 4$, the vertices $u$ and $w$ are adjacent.
Then $u$ is adjacent to 17 vertices of $X_1$ and to 4 vertices $w,x_1,x_2,x_3$,
hence, $|\Delta(u)|>19$, a contradiction. \wbull

\lmml{Case1}
Case $(1)$ is impossible.
\elmm
\proof We note that $|X_1|=|X_2|=16$ and $|X_3|=17$.

Let us first consider the case, when 
$X_1\subset \Delta(p)$, $X_2\subset \Delta(q)$. 
We may assume that $X_1\subset \Delta(u)$. 
Then $p\sim u$, $X_2\subset \Delta(w)$ and $q\sim w$.

Let $y_1$ be a vertex of $X_1$. Then $y_1$ is adjacent 
to 18 vertices of $X_1\cup \{x_1,p,u\}$, 
hence, $y_1$ is adjacent to a vertex of $X_2\cup X_3$.
Therefore, there are exactly 16 edges between $X_1$ and $X_2\cup X_3$.

Let $y_3$ be a vertex of $X_3$. Then $y_3$ is adjacent 
to 17 vertices of $X_3\cup \{x_3\}$, hence, $y_3$ is adjacent 
to a couple of vertices of $X_1\cup X_2$. 
Since $|\Delta(p,y_3)|\le 2$, the vertex $y_3$ has 
exactly one neighbor in $X_1$ and 
exactly one neighbor in $X_2$. 
This implies that there are 17 edges 
between $X_1$ and $X_3$, which is impossible.

We may now suppose $X_3\subset \Delta(p)$. 
Then $X_2\subset \Delta(q)$ and $p\not\sim q$. 

Suppose that $X_1\subset \Delta(u)$, $X_2\subset \Delta(w)$.
Then $q\sim w$, $w\not\sim u$ and $u$ is adjacent to a vertex of $X_3$.
Let $y_1$ be a vertex of $X_1$. 
Then $y_1$ is adjacent to 17 vertices of $X_1\cup \{x_1,u\}$,
hence, $y_1$ is adjacent to a couple of vertices of $X_2\cup X_3$.
Because $\Delta(y_1,x_2)$ contains $u$ and $\Delta(y_1,p)$ contains $x_1$, 
the vertex $y_1$ has exactly one neighbor in $X_2$ and
exactly one neighbor in $X_3$. 
Let $y_2$ be a vertex of $X_2$. 
Then $y_2$ is adjacent to 18 vertices of $X_2\cup \{x_2,w,q\}$,
hence, $y_2$ is adjacent to a vertex of $X_1\cup X_3$.
Since there are exactly 16 edges between $X_1$ and $X_2$,
the vertex $y_2$ has exactly one neighbor in $X_1$.
Let $a$ be a vertex of $\Gamma(x_1,w)\cap \Gamma_2(\infty)$.
By Lemma \ref{Claim3}, the vertex $a$ has a neighbor (say, $b$) 
in $X_1\cup X_2\cup \{u,p,x_2,q\}$.
If $b\in X_1\cup X_2$, then we have a contradiction with Lemma \ref{Claim3}.
If $b=u$ (or $b=x_2$), then $|\Gamma(u,w)|>3$ (or $|\Gamma(x_1,x_2)|>3$), 
a contradiction. Hence, $b\in \{p,q\}$.
Since $|\Gamma(x_1,w)\cap \Gamma_2(\infty)|=18$, we have 
$|\Gamma(x_1,q)|>3$ or $|\Gamma(p,w)|>3$, which is impossible.

At last, suppose that $X_1\subset \Delta(u,w)$ (and $u\sim w$), 
$X_2\subset \Delta(q)$ and $q$ is adjacent to a vertex of $X_1\cup X_3$.
Let $y_2$ be a vertex of $X_2$. 
Then $y_2$ is adjacent to 17 vertices of $X_2\cup \{x_2,q\}$,
hence, $y_2$ is adjacent to a couple of vertices of $X_1\cup X_3$.
Because $\Delta(y_2,x_3)$ contains $q$ and $\Delta(y_2,u)$ contains $x_2$, 
the vertex $y_2$ has exactly one neighbor in $X_3$ and exactly one neighbor, 
say $y_1$, in $X_1$. Now $y_1\not\sim x_2$ and $\Delta(y_1,x_2)=\{u,w,y_2\}$, 
a contradiction.
The lemma is proved.\wbull

\lmml{Case2}
Case $(2)$ is impossible.
\elmm
\proof We note that $|X_1|=|X_3|=17$ and $|X_2|=15$.

If $X_2\subset \Delta(u,w,p,q)$, then 
the vertices $u,w,p,q$ are mutually adjacent, which is impossible.
If $X_1\subset \Delta(u,w)$, then $u\sim w$ and 
the subgraph $\Delta(u)$ contains 17 vertices of $X_1$ and 
3 vertices $x_1,w,x_2$, hence, $|\Delta(u)|>19$, a contradiction.
So, we may assume that $X_1\subset \Delta(u)$,
$X_2\subset \Delta(w,q)$ and, hence, $w\sim q$.

Let us consider the case $X_3\subset \Delta(p)$.
Since $\Delta(p,q)=\{x_2,x_3\}$ and the subgraph $\Delta(q)$ 
contains $X_2\cup \{w,x_2,x_3\}$, 
the vertex $q$ is adjacent to a vertex of $X_1$ (say, $q'$).
Symmetrically, the vertex $w$ is adjacent to a vertex of $X_3$ (say, $w'$).
A vertex $y_1\in X_1-\{q'\}$ is adjacent to 18 vertices 
of $X_1\cup \{x_1,u\}$ and to a vertex of $X_2\cup X_3$. 
Similarly, a vertex $y_3\in X_3-\{w'\}$ is adjacent to 18 vertices 
of $X_3\cup \{x_3,p\}$ and to a vertex of $X_2\cup X_1$.
A vertex $y_2\in X_2$ is adjacent to 17 vertices of $X_2\cup \{x_2,w,q\}$
and to a couple of vertices of $X_1\cup X_3$.
Hence, there are the vertices, say, $a\in X_1$, $b\in X_3$ such that $a\sim b$ 
and each vertex $y_1\in X_1-\{a,q'\}$ ($y_3\in X_3-\{b,w'\}$, respectively)
has exactly one neighbor in $X_2$.
Further, the vertices $b,u,q$ induce a coclique of size 3 such that 
$\Delta(b,u)=\{a\}$, $\Delta(b,q)=\{x_3\}$ and $\Delta(u,q)=\{x_2,q'\}$ 
and this is Case (1), which is impossible.

At last, suppose that $X_2\subset \Delta(p)$, i.e., $X_2\subset \Delta(w,p,q)$. 
Then $p\sim q$ and $p\sim w$. 
Let $y_1$ be a vertex of $X_1$. Then $y_1$ is adjacent 
to 18 vertices of $X_1\cup \{x_1,u\}$, hence, $y_1$ is adjacent to 
a vertex of $X_2\cup X_3$. 
Let $y_3$ be a vertex of $X_3$. Then $y_3$ is adjacent 
to 17 vertices of $X_3\cup \{x_3\}$, hence, $y_3$ is adjacent 
to a couple of vertices of $X_1\cup X_2$.
Let $y_2$ be a vertex of $X_2$. Then $y_2$ is adjacent 
to 18 vertices of $X_2\cup \{x_2,w,p,q\}$, hence, $y_2$ 
is adjacent to one vertex of $X_1\cup X_3$. 
So, 
there are exactly 17 edges between $X_1$ and $X_2\cup X_3$,
exactly 15 edges between $X_2$ and $X_1\cup X_3$ and 
exactly 34 edges between $X_3$ and $X_1\cup X_2$, 
which is impossible. The lemma is proved.\wbull

\lmml{Case3}
Case $(3)$ is impossible.
\elmm
\proof We note that $|X_1|=|X_3|=17$, $|X_2|=16$
and $\Delta(u)=X_2\cup \{x_1,x_2,x_3\}$.
Hence, we have $X_1\subset \Delta(p)$ and $X_3\subset \Delta(q)$.
Let $y_1$ be a vertex of $X_1$. Then $y_1$ is adjacent 
to 18 vertices of $X_1\cup \{x_1,p\}$, hence, $y_1$ is adjacent to 
a vertex of $X_3\cup X_2$. 
Similarly, a vertex $y_3\in X_3$ is adjacent to 
a vertex of $X_1\cup X_2$.
Let $y_2$ be a vertex of $X_2$. Then $y_2$ is adjacent 
to 17 vertices of $X_2\cup \{x_2,u\}$, hence, $y_2$ 
is adjacent to two vertices of $X_1\cup X_3$. 
Therefore, each vertex $y_2\in X_2$ has exactly one 
neighbor in $X_1$ and one neighbor in $X_3$ and 
there is an edge $\{a,b\}$, where $a\in X_1$, $b\in X_3$.
Now the vertices $b,x_1,x_2$ induce a coclique of size 3
such that $\Delta(x_1,x_2)=\{p,u\}$, $\Delta(x_1,b)=\{a\}$ 
and $\Delta(x_2,b)=\{q\}$ and this is Case (1), a contradiction
with Lemma \ref{Case1}. 
This contradiction completes the proof of Theorem \ref{T}.


\end{document}